\documentclass{article}[12pt]

\usepackage{amssymb}
\usepackage[english]{babel}
\usepackage{amsfonts,epsfig}
\usepackage{epstopdf}
\addto\captionsenglish{}

\renewcommand{\i}{i}
\newcommand{\ptl}{\partial}

\newtheorem{definition}{Definition}[section]

\newtheorem{theorem}[definition]{Theorem}

\newtheorem{problem}{Problem}

\newcommand{\be}{\begin{equation}}
\newcommand{\ee}{\end{equation}}

\begin{document}

\title{Diffraction by an impedance strip II. Solving Riemann--Hilbert problems by OE--equation method}

\author{A. V. Shanin \and A. I. Korolkov}



\maketitle


\begin{abstract}
The current paper is the second part of a series of two papers
dedicated to 2D problem of  diffraction
of acoustic waves by a segment bearing impedance boundary conditions. In the first part some preliminary steps
were made, namely, the problem was reduced to two matrix Riemann--Hilbert problem. Here the Riemann--Hilbert
problems are solved with the help of a novel method of OE--equations.

Each Riemann--Hilbert problem is embedded into a family of similar problems with the same coefficient
and growth condition, but with some other cuts. The family is indexed by an artificial parameter. It is proven
that the dependence of the solution on this parameter can be described by a simple ordinary differential equation
(ODE1). The boundary conditions for this equation are known and the inverse problem of
reconstruction of the coefficient of ODE1 from the boundary conditions is formulated. This problem
is called the OE--equation. The OE--equation is solved by a simple numerical algorithm.
\end{abstract}

\section{Introduction}

This paper is the second part of a big work dedicated to
diffraction of a plane wave by a thin  infinite impedance strip.
In \cite{1} (which will be referred to as Part~I hereafter) some preliminary steps were made.
Namely, the diffraction
problem was formulated and symmetrized.
Functional problems
of the Wiener--Hopf class with entire functions
were introduced. Using the method of embedding formula these problems were reduced to two auxiliary problems.
Finally, two Riemann--Hilbert problems were formulated.

\begin{figure}[ht]
\centerline{\epsfig{file=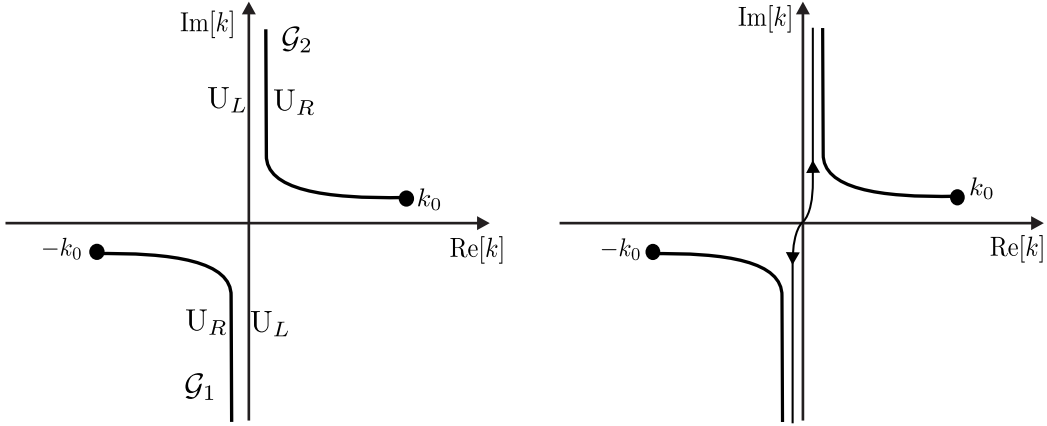}}
\caption{(left)Contours  ${\cal G}_{1,2}$  (right) Analytical continuation of
the square roots}
\label{fig00}
\end{figure}

\begin{figure}[ht]
\centerline{\epsfig{file=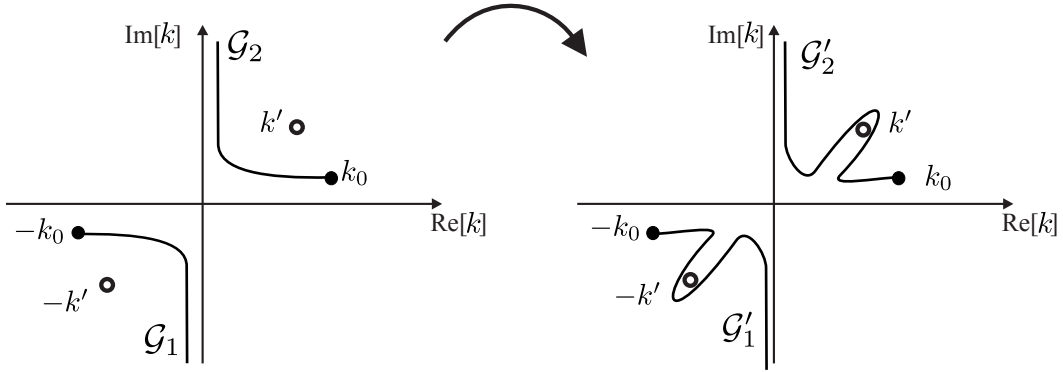,width = 14cm}}
\caption{Deformation of the cuts ${\cal G}_{1,2}$}
\label{fig05}
\end{figure}

The Rimann--Hilbert problems are formulated on the complex plane with cuts ${\cal G}'_{1,2}$. The cuts depend on
the impedance of the segment $\eta$. Due to energy absorption/conservation principle the impedance should obey
the condition ${\rm Im}[\eta] \le 0$. Then, if ${\rm Re}[\eta]  >0$ the contours ${\cal G}'_{1,2}$ coincide with
the undeformed contours ${\cal G}_{1,2}$ shown in Fig.~\ref{fig00} (left). These contours correspond to the trajectory of
the square root $\pm i\sqrt{k_0^2 - k^2}$ as $k$ takes real values (we remind that $k_0$ has a small positive imaginary
part). If ${\rm Re}[\eta]  \le 0$ the contours
${\cal G}'_{1,2}$ are obtained from ${\cal G}_{1,2}$ as the result of
deformation shown in Fig.~\ref{fig05}.
Points $\pm k'$ in the figure are zeros of $\eta - i \sqrt{k_0^2 - k^2}$:
$$
k' = \sqrt{k_0^2 + \eta^2}.
$$
The cuts ${\cal G}'_{1,2}$ are assumed to be symmetrical:
${\cal G}'_{2} = -{\cal G}'_{1}$.

The aim of the deformation shown in Fig.~\ref{fig05} is to make zeros of $\eta - i \sqrt{k_0^2 - k^2}$
not belonging to the plane cut along
${{\cal G}'_{1,2}}$.

For the antisymmetrical auxiliary problem the Riemann--Hilbert has form:

\begin{problem}
Find a matrix function
$${\rm U}(k) = \left(
\begin{array}{cc}
U^1_-(k) & U^1_+ (k) \\
U^2_-(k) & U^2_+ (k)
\end{array}
\right) $$
such that
\begin{itemize}
 \item it is regular on the complex plane cut along the lines ${{\cal G}'_{1,2}}$ (see Fig.~\ref{fig00}, left);

\item it obeys the following functional equations connecting the values on the shores of the cuts:
\begin{equation}
{\rm U}_R(k) = {\rm U}_L(k)\, {\rm M}_1(k), \qquad k \in {\cal G}'_1 ,
\label{eq0801}
\end{equation}
\begin{equation}
{\rm U}_R(k) = {\rm U}_L(k)\, {\rm M}_2(k), \qquad k \in {\cal G}'_2 ,
\label{eq0802}
\end{equation}
with coefficients
\begin{equation}
{\rm M}_1 (k) = \left(  \begin{array}{cc}
1 & 2 i \xi/(\eta - i \xi) \\
0 & (\eta + i \xi)/(\eta - i \xi)
\end{array} \right),
\label{eq0806}
\end{equation}
\begin{equation}
{\rm M}_2 (k) = \left(  \begin{array}{cc}
(\eta + i \xi)/(\eta - i \xi) & 0 \\
2 i \xi/(\eta - i \xi)      & 1
\end{array} \right).
\label{eq0805}
\end{equation}

\item it obeys the following growth restrictions:
\begin{equation}
U^{j}_+(k) = \delta_{j,2} (e^{-i \pi /2 } k)^{1/2} e^{i k a} + O(k^{-1/2} e^{i k a}),
\qquad {\rm Arg}[e^{- i \pi /2} k] \le \pi/2 ,
\label{eq0609c}
\end{equation}
\begin{equation}
U^{j}_-(k) = \delta_{j,1} (e^{i \pi /2 } k)^{1/2} e^{-i k a} + O(k^{-1/2} e^{-i k a}),
\qquad {\rm Arg}[e^{ i \pi /2} k] \le \pi/2 ,
\label{eq0610c}
\end{equation}
\begin{equation}
U^{j}_- (k) = i \, \delta_{j,1} ( e^{- i \pi /2 } k)^{1/2} e^{- i k a} + O(k^{-1/2} e^{-i k a}),
\qquad {\rm Arg}[e^{- i \pi /2} k] \le \pi/2 ,
\label{eq0807}
\end{equation}
\begin{equation}
U^{j}_+ (k)= i \, \delta_{j,2} ( e^{  i \pi /2 } k)^{1/2} e^{i k a} + O(k^{-1/2} e^{i k a}),
\qquad {\rm Arg}[e^{ i \pi /2} k] \le \pi/2 .
\label{eq0808}
\end{equation}

\item functions $U^j_{\pm}$ grow no faster than a constant near the points $\pm k_0$.
\end{itemize}
\label{WHH}
\end{problem}

Notations ${\rm U}_{L,R}$ correspond to the values of ${\rm U}$ taken on the left and right shores of the cuts
(see Fig.~\ref{fig00});
\[
\xi = \xi(k) \equiv \sqrt{k_0^2 - k^2}.
\]
Square root is equal to $k_0$ at the point $k=0$ and then continued to  ${{\cal G}'_{1,2}}$,
along the contours shown in Fig.~\ref{fig00}, right.

For the symmetrical case the Riemann--Hilbert problem has form:
\begin{problem}
Find a matrix function
$${\rm V}(k) = \left(
\begin{array}{cc}
V^{1}_-(k) & V^{1}_+ (k) \\
V^{2}_-(k) & V^{2}_+ (k)
\end{array}
\right) $$
such that
\begin{itemize}
\item it is regular on the plane cut along the lines ${{\cal G}'_{1,2}}$;

\item it obeys functional equations
\begin{equation}
{\rm  V}_R(k) = {\rm  V}_L {\rm N}_1(k), \qquad k \in {\cal G}'_1 ,
\label{eq0801sym}
\end{equation}
\begin{equation}
{\rm  V}_R(k) = {\rm  V}_L {\rm N}_2(k), \qquad k \in {\cal G}'_2
\label{eq0802sym}
\end{equation}
with coefficients
 \begin{equation}
{\rm  N}_1 (k) = \left(  \begin{array}{cc}
1 & -2 \eta/(\eta - i \xi) \\
0 & (\eta + i \xi)/(i \xi - \eta  )
\end{array} \right) ,
\label{eq0806sym}
\end{equation}
\begin{equation}
{\rm N}_2 (k) = \left(  \begin{array}{cc}
(\eta + i \xi)/(i \xi - \eta) & 0 \\
-2 \eta/(\eta - i \xi)      & 1
\end{array} \right) ,
\label{eq0805sym}
\end{equation}
on the cuts;

\item it obeys the following growth restrictions:
\begin{equation}
V^{j}_+(k) = \delta_{j,2}  e^{i k a} + O(k^{-1}e^{i k a}),
\qquad {\rm Arg}[e^{- i \pi /2} k] \le \pi/2 ,
\label{eq0609s}
\end{equation}
\begin{equation}
V^{j}_-(k) = \delta_{j,1}  e^{-i k a} + O(k^{-1} e^{-i k a}),
\qquad {\rm Arg}[e^{ i \pi /2} k] \le \pi/2 ,
\label{eq0610s}
\end{equation}
\begin{equation}
 V^{j}_- (k)=  \, \delta_{j,1}e^{- i k a} + O(k^{-1}e^{-i k a}),
\qquad {\rm Arg}[e^{- i \pi /2} k] \le \pi/2  ,
\label{eq0807sym}
\end{equation}
\begin{equation}
 V^{j}_+ (k)=  \, \delta_{j,2} e^{i k a} + O(k^{-1} e^{i k a}),
\qquad {\rm Arg}[e^{ i \pi /2} k] \le \pi/2  .
\label{eq0808sym}
\end{equation}

\item functions $ V^{j}_{\pm}$ grow no faster than $(\sqrt{k_0 \mp k})^{-1/2}$ near the points $\pm k_0$.
\end{itemize}
\label{WHH_sym}
\end{problem}

If we manage to find a solution of Problem~\ref{WHH}, we can recover a antisymmetrical part of
the solution of original problem using following procedure. First, functions $\tilde U_0^1(k), \tilde U_0^2(k)$ are calculated:
\begin{equation}
\tilde U_0^j = -(\eta - i\sqrt{k_0^2 - k^2})^{-1}(U^j_- + U^j_+),\qquad j =1,2,
\label{u0tildej}
\end{equation}
Then function $\tilde U^0(k,k_*)$ is found by the embedding formula:
 \begin{equation}
\tilde U_{0}(k,k_*) =\frac{\xi(k_*)}{k-k_*}\left(\tilde U_0^{1}(k_*)\tilde U^{2}_0(k)  - \tilde U_0^{1}(k)\tilde U^{2}_0(k_*)\right).
\label{embedding_a}
\end{equation}
$$
k_* = k_0 \cos\theta^{\rm in}.
$$
Finally, the antisymmetrical part of the directivity is found:
\begin{equation}
S^{\rm a} (\theta,\theta^{\rm in}) = -e^{- i \pi /4} k_0 \sin \theta \,  \tilde U_0 (- k_0 \cos\theta).
\label{eq0617}
\end{equation}

For the symmetrical case (Problem~\ref{WHH_sym}) the following formulae are used:
\begin{equation}
\tilde V^j_0 =  \frac{\xi(k)}{i(\eta - i \xi(k) )} (V^j_- + V^j_+), \qquad j=1,2,
\end{equation}
\begin{equation}
\tilde V_{0}(k,k_*) =\frac{i\eta}{(k-k_*)}\left(\tilde V^{2}_0(k_*)\tilde V_0^{1}(k) - \tilde V^{2}_0(k)\tilde V_0^{1}(k_*)  \right),
\label{embedding_s}
\end{equation}
\begin{equation}
S^{\rm s} (\theta, \theta^{\rm in}) =  e^{- i \pi /4}  \tilde V_0 (- k_0 \cos(\theta)).
\label{eq0617a}
\end{equation}
The directivity related to the initial problem is a sum of the antisymmetrical and symmetrical part:
\begin{equation}
S(\theta, \theta^{\rm in}) = S^{\rm a}(\theta, \theta^{\rm in}) +
S^{\rm s}(\theta, \theta^{\rm in}).
\end{equation}

In the present paper we solve Problem~\ref{WHH} and Problem~\ref{WHH_sym}.
We use for this the method of OE--equation proposed  recently.
The plan of the research is as follows.
First, a family of Riemann--Hilbert problems indexed by an artificial parameter $b$ is formulated. Then, an ordinary differential equation with respect to $b$ (ODE1) is introduced. This equation is supplemented with initial conditions. An OE--equation (an equation for the coefficients of ODE1) is formulated. This equation is solved numerically; the results  are compared with solutions obtained by the integral equation method.


\section{A family of Riemann--Hilbert problems}

\subsection{One more preliminary step for the antisymmetrical case}

A crucial step of the method is introducing of the family of Riemann-Hilbert problems to which
Problems~\ref{WHH} and \ref{WHH_sym} belong as an element. Before we introduce such a family
it is necessary to reformulate the Riemann--Hilbert problems (Problem~1 and~2)
in such a way
that the connection matrices ${\rm M}_{1,2}(k)$ and ${\rm N}_{1,2}(k)$ have eigenvalues
tending to 1 as $|k|\to \infty$. One can see that matrices
${\rm N}_{1,2}$ satisfy this condition (so no reformulation is needed),
while matrices ${\rm M}_{1,2}(k)$ have one eigenvalue  tending to 1, and the other tending to $-1$.
To reformulate the antisymmetrical problem make the following variable change:
\begin{equation}
{\rm \hat U} \equiv
\left( \begin{array}{cc} \hat U^1_- & \hat U^1_+ \\
 \hat U^2_- & \hat U^2_+ \end{array} \right) =
\left( \begin{array}{cc} U^{1}_- & U^{1}_+ \\
 U^{2}_- & U^{2}_+ \end{array} \right)
\left( \begin{array}{cc} e^{i \pi/4} (k_0 - k)^{-1/2} & 0 \\
 0 & e^{i \pi/4} (k_0 + k)^{-1/2} \end{array} \right) .
\label{eq0809}
\end{equation}

The growth restrictions for the new functions become as follows:
\begin{equation}
\hat U^{j}_+(k) = \delta_{j,2} e^{i k a} + O(k^{-1} e^{i k a}),
\qquad {\rm Arg}[e^{ -i \pi /2} k] \le \pi/2 ,
\label{eq0810}
\end{equation}
\begin{equation}
\hat U^{j}_-(k) = \delta_{j,1} e^{-i k a} + O(k^{-1} e^{-i k a}),
\qquad {\rm Arg}[e^{ i \pi /2} k] \le \pi/2 ,
\label{eq0811}
\end{equation}
\begin{equation}
\hat U^{j}_-(k) = - \delta_{j,1}  e^{- i k a} + O(k^{-1} e^{-i k a}),
\qquad {\rm Arg}[e^{ -i \pi /2} k] \le \pi/2 ,
\label{eq0812}
\end{equation}
\begin{equation}
\hat U^{j}_+(k) =  - \delta_{j,2}  e^{i k a} + O(k^{-1} e^{i k a}),
\qquad {\rm Arg}[e^{ i \pi /2} k] \le \pi/2  .
\label{eq0813}
\end{equation}

The connection formulae for ${\rm \hat U}$ on the cuts become as follows:
\begin{equation}
{\rm \hat U}_R(k) = {\rm \hat U}_L \tilde {\rm M}_1(k), \qquad k \in {\cal G}'_1 ,
\label{eq0814}
\end{equation}
\begin{equation}
{\rm \hat U}_R(k) = {\rm \hat U}_L \tilde {\rm M}_2(k), \qquad k \in {\cal G}'_2 ,
\label{eq0815}
\end{equation}
where
\begin{equation}
\tilde {\rm M}_2 (k) = \left(  \begin{array}{cc}
(i \xi + \eta)/( i \xi- \eta ) & 0 \\
2 i (k_0 + k) /(i \xi- \eta)      & 1
\end{array} \right) ,
\label{eq0816}
\end{equation}
\begin{equation}
\tilde {\rm M}_1 (k) = \left(  \begin{array}{cc}
1 & 2 i (k_0 - k)/( i \xi- \eta)   \\
0 & (i \xi + \eta )/( i \xi - \eta )
\end{array} \right) .
\label{eq0817}
\end{equation}

We can formulate now a functional problem for ${\rm \hat U}$, which replaces Problem~\ref{WHH}:

\begin{problem}
Find a matrix function ${\rm \hat U}(k)$ of elements (\ref{eq0809}) such that

\begin{itemize}
\item it is regular and has no zeros of determinant on the plane cut along the lines ${\cal G}'_{1,2}$;

\item it obeys functional equations (\ref{eq0814}), (\ref{eq0815}) with coefficients
(\ref{eq0816}), (\ref{eq0817}) on the cuts;

\item it obeys growth restrictions (\ref{eq0812}), (\ref{eq0813}), (\ref{eq0814}), (\ref{eq0815});

\item components of ${\rm \hat U}$ grow no faster than $(k_0\mp k)^{-1/2}$  near the points $\pm k_0$.
\end{itemize}
\label{WHH_V}
\end{problem}

\subsection{A family of Riemann--Hilbert problems in the antisymmetrical case}

Consider the antisymmetrical case, i.\ e.\
Problem~\ref{WHH_V}.
Represent contours ${\cal G}_{1,2}'$ as ${\cal G}_2' = \gamma+ k_0$, ${\cal G}_1' =- \gamma-k_0$, where
$\gamma$ is a contour going from $i \infty$ to $0$. Here $+k_0$ or $-k_0$ means a shift of the contour.

Let $\gamma(b)$, $b \in \gamma$, be a contour going from $i \infty$ to $b$ along $\gamma$. I.\ e.\ $\gamma(b)$ is a part of $\gamma$. Let be
\[
{\cal G}_1'(b) =- \gamma(b) -k_0 ,
\qquad
{\cal G}_2'(b) =  \gamma(b)+k_0.
\]

The family of the Riemann--Hilbert problems is built based of Problem~\ref{WHH_V}. The key step
is to replace the contours ${\cal G}'_{1,2}$ with ${\cal G}'_{1,2} (b)$. The growth conditions
at infinity and the connection matrices
(\ref{eq0816}), (\ref{eq0817}) remain the same as for Problem~\ref{WHH_V}, while the growth restrictions
at the ends of the contours should be changed (since the ends of the contours change from $\pm k_0$ to
$\pm (k_0 + b)$). To formulate these restrictions we should study behavior of a solution of the Riemann--Hilbert
problem near the end of one of the contours. For example, consider contour ${\cal G}_1'(b)$. Let equation
(\ref{eq0815}) with coefficient (\ref{eq0817}) be fulfilled on its shores. Then, obviously, at the vicinity of the
end point $b+k_0$ the solution has form
\begin{equation}
{\rm \hat U} = {\rm T}(k) \, {\rm H} \,
\exp \left\{ \frac{\log (k - (k_0 + b))}{2 \pi i}
\left( \begin{array}{cc}
\log (m_1(b))  &  0 \\
 0 & \log (m_2(b))
\end{array} \right) \right\} \, {\rm H}^{-1}
\label{eq0906}
\end{equation}
where ${\rm T} (k)$ is an arbitrary matrix analytical near $k_0 + b$,
\begin{equation}
m_1(b) \equiv m(b) = \frac{i\sqrt{k_0^2 - (k_0 + b)^2} + \eta}{i \sqrt{k_0^2 - (k_0 + b)^2} - \eta},
\qquad
b \in \gamma .
\label{eq0901}
\end{equation}
and
\[
m_2 (b) = 1
\]
are eigenvalues of $\tilde {\rm M}_2(k_0 + b)$,
\begin{equation}
{\rm H} = \left( \begin{array}{cc}
1  & 0 \\
\alpha & 1
\end{array} \right) ,
\qquad
\alpha =
\frac{- i (k_0 + k)}{\eta}.
\label{eq0908}
\end{equation}
is the matrix of eigenvectors of
$\tilde {\rm M}_2(k_0 + b)$.
The branch of the square root at $k = b + k_0$ is chosen according to the explanation above.

An appropriate choice of the logarithms in (\ref{eq0906}) determines the growth restrictions near $k_0 + b$.
Choose $\log (m_2 (b)) =0$ (this corresponds to a regular component of the solution). Then, consider the function $m(b)$.
Obviously $m (0) = -1$,  $m(i \infty) = 1$.

Introduce an important value
\begin{equation}
{\rm Idx} = \left. \log(m (b)) \right| _0^{i \infty}  ,
\label{eq0902}
\end{equation}
which will be called {\em the index\/} of the Riemann--Hilbert problem discussed here.
The notation above denotes the continuous change of the logarithm value along the contour
${\cal G}'_2$.

Obviously, ${\rm Idx} = \pi i + 2 \pi i \mu$ for some integer $\mu$.
It is not difficult to show (see Appendix~A) that under the restriction ${\rm Im}[\eta] < 0$
\begin{equation}
{\rm Idx} = \pi i .
\label{eq0903}
\end{equation}

Define also the value
\begin{equation}
\lambda(b) = \frac{\log(m(b))}{2 \pi i} .
\label{eq0904}
\end{equation}
This function should be continuous on ${\cal G}'_2$, and besides
\begin{equation}
\lambda(i \infty) = 0.
\label{eq0905}
\end{equation}
According to (\ref{eq0903}),
 $\lambda(0) = -1/2$.

Introduce a family of $2 \times 2$ matrix
functions ${\rm \hat U}(b,k)$ such that for any fixed $b$ the function ${\rm \hat U}(b,k)$ taken as the
function of $k$ is a solution of the following functional problem:
\begin{problem}
Find a matrix function ${\rm \hat U}(b,k)$ with elements denoted as (\ref{eq0809}) such that
\begin{itemize}
\item it is regular and has no zeros of determinant on the plane cut along the lines ${\cal G}_{1,2}'(b)$;

\item it obeys functional equations (\ref{eq0814}), (\ref{eq0815}) with coefficients
(\ref{eq0816}), (\ref{eq0817}) on the cuts ${\cal G}_{1,2}'(b)$;

\item it obeys growth restrictions (\ref{eq0810}), (\ref{eq0811}), (\ref{eq0812}), (\ref{eq0813})
at infinity;

\item near $-k_0 - b$ components $\hat U^j_-$ behave as regular functions, components
   $\hat U^j_+ $ behave as
$(k-(k_0 +b))^{\lambda(b)} \Psi_1(k+(k_0 +b)) + \Psi_2(k+(k_0 +b))$,
where $\Psi_1$ and $\Psi_2$ are some functions regular near zero;

\item near $k_0 + b$ components $\hat U^j_+$ behave as regular functions, components $\hat U^j_-$
behave as
$(k+(k_0 +b))^{\lambda(b)} \Psi_3(k-(k_0 + b)) + \Psi_4(k - (k_0 + b))$,
where $\Psi_3$ and $\Psi_4$ are some functions regular near zero.
\end{itemize}
\label{family}
\end{problem}

The definition of ${\rm \hat U}(b, k)$ is mathematically correct, since uniqueness of ${\rm \hat U}(b,k)$
can be proven for each $b$. The proof is based on the determinant technique introduced in Part~I.

Problem~\ref{WHH_V} and Problem~\ref{family} are connected via the relation
\begin{equation}
{\rm \hat U}(k) = {\rm \hat U}(0, k).
\label{eq0907}
\end{equation}


\subsection{A family of Riemann--Hilbert problems in the symmetrical case}

In the symmetrical case introduce  a family of functions ${\rm V}(b,k)$ such that for any fixed $b$ the function ${\rm V}(b,k)$ taken as the
function of $k$ is a solution of the following functional problem:
\begin{problem}
Find a matrix function ${\rm V}(b,k)$ such that
\begin{itemize}

\item it is regular and has no zeros of determinant on the plane cut along the lines ${\cal G}_{1,2}'(b)$;

\item it obeys functional equations (\ref{eq0801sym}), (\ref{eq0802sym}) with coefficients
(\ref{eq0806sym}), (\ref{eq0805sym}) on the cuts ${\cal G}_{1,2}'(b)$;

\item it obeys growth restrictions (\ref{eq0609s}), (\ref{eq0610s}), (\ref{eq0807sym}), (\ref{eq0808sym})
at infinity;

\item near $-k_0 - b$ components $V^{j}_-$ behave as regular functions, components
   $V^{j}_+ $ behave as
$(k+(k_0 +b))^{\lambda(b)} \tilde \Psi_1(k+(k_0 +b)) + \tilde \Psi_2(k+(k_0 +b))$,
where $\tilde \Psi_1$ and $\tilde \Psi_2$ are some functions regular near zero;

\item near $k_0 + b$ components $V^{j}_+$ behave as regular functions, components $V^{j}_-$
behave as
$(k-(k_0 +b))^{\lambda(b)} \tilde \Psi_3(k-(k_0 + b)) + \tilde \Psi_4(k - (k_0 + b))$,
where $\tilde \Psi_3$ and $\tilde \Psi_4$ are some functions regular near zero.
\end{itemize}
\label{familysym}
\end{problem}
Problem~\ref{WHH_sym} and Problem~\ref{familysym} are connected via the relation
\begin{equation}
{\rm V}(k) = {\rm V}(0, k).
\label{eq0907sym}
\end{equation}


\section{Derivation of ODE1}
\subsection{ODE1 for the antisymmetrical problem}
We are looking for an ordinary differential equation (ODE1) in the form
\begin{equation}
\frac{\ptl }{\ptl b} {\rm \hat U}(b,k) =
{\rm R}(b, k) \, {\rm \hat U}(b,k),
\qquad
b \in \gamma,
\label{eq1001}
\end{equation}
where ${\rm R}(b, k)$ is the coefficient of the equation. Indeed, this equation is useful only if the
coefficient ${\rm R}$ has structure simpler than that of ${\rm \hat U}$. The form of the coefficient is given by the
following theorem.

\begin{theorem}
Function ${\rm \hat U}(b,k)$, which is a solution of a family of functional problems introduced as
Problem~\ref{family}, obeys equation (\ref{eq1001}) with the coefficient
\begin{equation}
{\rm R}(b,k) = \frac{{\rm r}(b)}{k-(k_0 + b)} -
\frac{{\rm r}^*(b)}{k+(k_0 + b)},
\label{eq1002}
\end{equation}
where ${\rm r}(b)$, $b \in \gamma$ is a $2 \times 2$ matrix function (not depending on $k$);
${\rm r}^*$ is connected with  ${\rm r}$ via the relation:
\[
{\rm r} = \left( \begin{array}{cc}
r_{1,1} & r_{1,2} \\
r_{2,1} & r_{2,2}
\end{array} \right),
\qquad
{\rm r}^* = \left( \begin{array}{cc}
r_{2,2} & r_{2,1} \\
r_{1,2} & r_{1,1}
\end{array} \right),
\]
i.\ e.\ to obtain ${\rm r}^*$ one has to interchange first the rows and then the columns of~${\rm r}$.
\label{ODE1}
\end{theorem}

\noindent
{\bf Proof}
Construct the coefficient of ODE~1 as follows:
\begin{equation}
{\rm R}(b, k) = \frac{\ptl {\rm \hat U}(b,k)}{\ptl b} {\rm \hat U}^{-1}(b,k).
\label{eq1003}
\end{equation}
Consider this combination for fixed $b$ as the function of $k$. According to Problem~\ref{family},
${\rm R}(b,k)$ has no singularities on the complex plane $k$ cut along the contours
${\cal G}_{1,2}' (b)$. Moreover, since the functions $\tilde {\rm M}_{1,2}$ do not depend on $b$, the
 values of ${\rm R}$ on the left and right shores of ${\cal G}_{1,2}'(b)$ are equal:
\[
\left( \frac{\ptl {\rm \hat U}(b,k)}{\ptl b} {\rm \hat U}^{-1}(b,k) \right)_R =
\left( \frac{\ptl {\rm \hat U}(b,k)}{\ptl b}\right)_R  \left( {\rm \hat U}^{-1}(b,k) \right)_R =
\]
\[
\frac{\ptl ({\rm \hat U}(b,k))_L \tilde {\rm M}_j(k) }{\ptl b} \tilde {\rm M}_j(k)^{-1} ({\rm \hat U}^{-1}(b,k))_L =
\left(\frac{\ptl {\rm \hat U}(b,k)}{\ptl b} \right)_L ({\rm \hat U}^{-1}(b,k))_L.
\]
Thus, function ${\rm R}$ is single-valued on the plane $k$. The only singularities it can have
are the ends of the contour ${\cal G}_{1,2}'(b)$, i.\ e.\ the points $k = \pm (k_0 + b)$.
Consider the vicinity of the point $k = k_0 + b$. To study the singularity of ${\rm R}$ at this point,
represent the solution ${\rm \hat U}$ in the form
\begin{equation}
{\rm \hat U}(b,k) = {\rm T}(b, k) {\rm H}(k) \exp \left\{
\log(k-(k_0 +b))  \left( \begin{array}{cc}
\lambda(k) & 0 \\
0     &  0
\end{array} \right)
 \right\} {\rm H}^{-1}(k) ,
\label{eq1004}
\end{equation}
where ${\rm T}(b,k)$ is a regular function in the area considered, ${\rm H}$ is the same function as
(\ref{eq0908}). Substituting (\ref{eq1004}) into (\ref{eq1003}),
obtain that the coefficient ${\rm R}$ can only have a simple pole at $k= k_0 + b$. Let ${\rm r}(b)$ be the
residue of ${\rm R}(b, k)$ at $k = k_0 + b$.

Similar consideration can be performed with respect to the point $k = - k_0 - b$.
Consider the geometrical symmetry of the problem, namely
$x \to -x$. This symmetry transforms the matrix of solutions as follows:
\begin{equation}
\left( \begin{array}{cc}
\hat U^1_-(k) & \hat U^1_+(k) \\
\hat U^2_-(k) & \hat U^2_+(k) \\
\end{array}\right)
=
\left( \begin{array}{cc}
\hat U^2_+(-k) & \hat U^2_-(-k) \\
\hat U^1_+(-k) & \hat U^1_-(-k) \\
\end{array}\right).
\label{eq1005}
\end{equation}
Due to this transformation, the coefficient ${\rm R}(b, k)$ has form of (\ref{eq1002}).
$\square$


\subsection{ODE1 for the symmetrical problem}

Similarly to the antisymmetrical case one can prove the following theorem.

\begin{theorem}
Function ${\rm V}(b,k)$, which is a solution of a family of functional problems introduced as
Problem~\ref{familysym}, obeys equation
\begin{equation}
\frac{\ptl }{\ptl b} {\rm V}(b,k) =
{\rm L}(b, k) \, {\rm V}(b,k),
\qquad
b \in \gamma,
\label{eq1001s}
\end{equation}
with the coefficient
\begin{equation}
{\rm L}(b,k) = \frac{{\rm l}(b)}{k-(k_0 + b)} -
\frac{{\rm l}^*(b)}{k+(k_0 + b)},
\label{eq1002s}
\end{equation}
where ${\rm} l(b)$, $b \in \gamma$ is a $2 \times 2$ matrix function (nor depending on $k$);
operator $\cdot^*$ is as introduced above
\label{ODE1s}
\end{theorem}


\subsection{Initial condition for ODE1}

\begin{theorem}
Initial conditions for ODE1 (\ref{eq1001}) and (\ref{eq1001s}) are as follows:
\begin{equation}
\lim_{b \to i \infty}{\rm \hat U}(b , k) =\lim_{b \to i \infty}{\rm V}(b , k) =
 \Pi(k),
\label{eq1006}
\end{equation}
\[
\Pi(k) =
\left( \begin{array}{cc}
\exp\{ - i a k  \} & 0 \\
0 & \exp\{  i a k  \}
\end{array} \right).
\]
\label{InitialConditions}
\end{theorem}

\noindent
{\bf Proof}
Consider the antisymmetrical case, i.\ e.\
consider Problem~\ref{family} for some large positive imaginary $b$. The functional problem can be
reduced to a system of integral equations as follows. Introduce the matrix
\begin{equation}
\left( \begin{array}{cc}
v^1_- & v^1_+ \\
v^2_- & v^2_+
\end{array} \right)
=
\left( \begin{array}{cc}
\hat U^1_- & \hat U^1_+ \\
\hat U^2_- & \hat U^2_+
\end{array} \right)
\Pi(k)
\label{eq1007}
\end{equation}
Then introduce the functions $\psi_+^j (k)$ $k \in {\cal G}'_1(b)$ and $\psi_-^j (k)$ $k \in {\cal G}'_2(b)$,
such that
\begin{equation}
v_-^j (k)= \delta_{j, 1} + \int_{{\cal G}_2'(b)} \frac{\psi_-^j (\tau)}{k - \tau} d\tau ,
\label{eq1008}
\end{equation}
\begin{equation}
v_+^j (k)= \delta_{j, 2} + \int_{{\cal G}_1'(b)} \frac{\psi_+^j (\tau)}{k - \tau} d\tau .
\label{eq1009}
\end{equation}
Assume that contours ${\cal G}_{1,2}(b)$ go from $\mp i \infty$ to $\mp (k_0 + b)$.
Note that for $k \in {\cal G}_2'(b)$
\begin{equation}
(v_-^j (k))_{L,R}= \delta_{j, 1} \pm \pi i \psi_-^j (k) +
\int_{{\cal G}_2'(b)} \frac{\psi_-^j (\tau)}{k - \tau} d\tau ,
\label{eq1008a}
\end{equation}
where the integral has sense of the main value.

According to the functional equation (\ref{eq0802}), the following equation is valid:
\[
-2\pi i \psi_-^j (k) =  e^{2 i a k} \tilde m_{2,1}(k) v^j_+ (k) +
\]
\begin{equation}
(\tilde m_{1,1}(k) -1)
\left( \int_{{\cal G}_2'(b)} \frac{\psi_-^j (\tau) d\tau}{k - \tau} + \pi i \, \psi_-^j (k)+ \delta_{1,j} \right),
\qquad
k \in {\cal G}_2'(b)
\label{eq1009a}
\end{equation}
According to the geometrical symmetry,
\begin{equation}
v_+^j (k) = \delta_{j,2} + \int_{{\cal G}_2'(b)} \frac{\psi^{3-j}_-(\tau)d\tau}{k+\tau} ,
\qquad
k \in {\cal G}_2'.
\label{eq1010}
\end{equation}
Here $\tilde m_{1,1}$ and $\tilde m_{2,1}$ are corresponding elements of matrix $\tilde {\rm M}_2$.
Note that for large ${\rm Im}[k]$, $k \in {\cal G}'_{2}(b)$, the values of $(\tilde m_{1,1} (k)-1)$ are close to~0
(this is the reason for variable change from ${\rm U}$ to ${\rm \hat U}$). Also, under the same condition
$e^{2 i a k} \tilde m_{2,1} (k)$ is close to zero.

If ${\rm Im}[b]$ is large enough, the system (\ref{eq1009a}) can be solved by iterations.
For large ${\rm Im}[b]$ only the zero-order approximation can be left, i.e. $\psi_\pm^j$ can be set to zero. This
corresponds to (\ref{eq1006}).

In the symmetrical case the proof is similar.  $\square$


\section{OE--equation}

\subsection{OE--notation}

Introduce the following notation.
Consider matrix ODE
\begin{equation}
\frac{\ptl }{\ptl \tau} {\rm X}(\tau) = {\rm K}(\tau)\, {\rm X}(\tau)
\label{eq1101}
\end{equation}
taken on a contour $h$ with starting point $\tau_1$ and ending point $\tau_2$. Let the
initial condition have form
\[
{\rm X}(\tau_1) = {\rm I}.
\]
By definition,
\begin{equation}
{\rm OE}_h\left[ {\rm K} (\tau)\, d\tau  \right]
\equiv
{\rm X} (\tau_2).
\label{eq1102}
\end{equation}

The following properties are obvious:
\begin{itemize}
\item  If $h'$ is the contour $h$ passed in the opposite direction, then
\begin{equation}
{\rm OE}_{h'}\left[ {\rm K} (\tau)\, d\tau  \right]
=\left(
{\rm OE}_{h}\left[ {\rm K} (\tau)\, d\tau  \right]
\right)^{-1}
\label{eq1103}
\end{equation}

\item If $h$ is a concatenation of $h_1$ and $h_2$ ($h_1$ is the first) then
\begin{equation}
{\rm OE}_{h}\left[ {\rm K} (\tau)\, d\tau  \right]
=
{\rm OE}_{h_2}\left[ {\rm K} (\tau)\, d\tau  \right]
{\rm OE}_{h_1}\left[ {\rm K} (\tau)\, d\tau  \right]
\label{eq1104}
\end{equation}
\end{itemize}


\subsection{Derivation of the OE--equation in the antisymmetrical case}

According to Theorems~\ref{ODE1} and~\ref{InitialConditions}, the solution of Problem~\ref{WHH_V}
can be written as
\begin{equation}
\hat {\rm U}(k) = {\rm OE}_{\gamma} \left[ \left(
\frac{{\rm r}(b)}{k - (k_0 + b)}
-
\frac{{\rm r^*}(b)}{k + (k_0 + b)}
\right) db \right] \Pi(k),
\label{eq1201}
\end{equation}

We remind that contour $\gamma$ goes from $i \infty$ to $0$.

A detailed study based on the continuation of matrices $\tilde {\rm M}_{1,2}$  near the
contours $\tilde {\cal G}_{1,2}'$ (see \cite{Shanin}) shows that the coefficient ${\rm R}$
is analytical with respect to the variable $b$ in a narrow strip surrounding contour $\gamma$. Thus,
the contour can be slightly deformed without changing the result, provided that the starting and the ending points of the contour remain the same.

\begin{figure}[ht]
\centerline{\epsfig{file=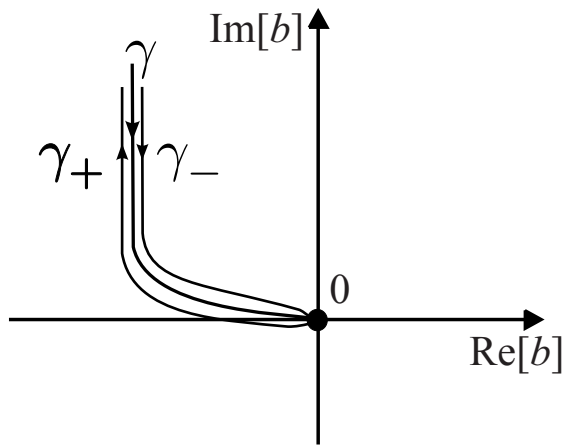}}
\caption{Contours $\gamma_+$ and $\gamma_-$. For simplicity the case ${\rm Re}[\eta] > 0$
is shown, i.\ e.\ the cuts ${\cal G}_{1,2}$ remain undeformed  }
\label{fig06}
\end{figure}

Draw contours $\gamma_+$ and $\gamma_-$ as it is shown in Fig.~\ref{fig06}. These contours are needed to
calculate the values ${\rm \hat U}_R(k)$ and ${\rm \hat U}_L(k)$, $k \in {\cal G}_2'$ without allowing
singulatities of the coefficient of ODE~1.
Namely,
\begin{equation}
\hat {\rm U}_R(k) = {\rm OE}_{\gamma_+} \left[
{\rm R}(b, k)
\, db \right] \Pi(k),
\label{eq1202}
\end{equation}
\begin{equation}
\hat {\rm U}_L(k) = \left( {\rm OE}_{\gamma_-} \left[
{\rm R}(b, k)
\, db \right] \right)^{-1} \Pi(k),
\label{eq1203}
\end{equation}

Define contour $\gamma_+ \circ \gamma_-$
as a concatenation of contours $\gamma_+$ and $\gamma_-$ ($\gamma_+$ is the first).
According to functional equation (\ref{eq0815}), the following relation is valid:
\begin{equation}
\Pi^{-1} (k) \, {\rm OE}_{\gamma_+ \circ \gamma_-} \left[ \left(
\frac{{\rm r}(b)}{k - (k_0 + b)}
-
\frac{{\rm r^*}(b)}{k + (k_0 + b)}
\right) db \right] \Pi(k) = \tilde {\rm M}_2 (k).
\label{eq1204}
\end{equation}
This is the OE--equation for the considered problem. Due to the geometrical symmetry, equation (\ref{eq0814}) will be also valid.

Formulate the problem for the OE--equation, to which the antisymmetrical problem of diffraction by
an impedance strip becomes reduced.

\begin{problem}
Find function ${\rm r}(b)$ for $b \in \gamma$ analytical in a narrow strip surrounding
$\gamma$ such that equation (\ref{eq1204}) is valid for each $k \in {\cal G}_2'$.
\label{OE}
\end{problem}

\subsection{OE--equation in the symmetrical case}

The following problem should be solved in the symmetrical case.

\begin{problem}
Find function ${\rm l}(b)$ analytical in a narrow strip surrounding
$\gamma$ such that equation
\begin{equation}
 \Pi^{-1} (k) \, {\rm OE}_{\gamma_+ + \gamma_-} \left[ \left(
\frac{{\rm  l}(b)}{k - (k_0 + b)}
-
\frac{{\rm  l^*}(b)}{k + (k_0 + b)}
\right) db \right]  \Pi(k) = {\rm  N}_2 (k)
\label{eq1204sym}
\end{equation}
be valid for each $k \in {\cal G}_2' $.
\label{OEsym}
\end{problem}


\section{Numerical results}


\subsection{Antisymmetrical case}

Solving of the diffraction problem by means of the proposed technique comprises the following steps:

\begin{itemize}

\item
Contour $\gamma$ (see Fig.~\ref{fig06}) is discretized. Problem~\ref{OE} is solved numerically (the procedure
is described below). As the result,
the coefficient ${\rm r}(b)$ becomes known in a set of points $b_n$ covering contour~$\gamma$ densely .

\item
Points of interest are selected in the $k$-plane. A good choice is the set of points $\kappa_m$ densely covering the
segment $(-k_0 , k_0)$, since such points enable one to construct the directivity of the field.
For these points the values $\tilde {\rm U}(\kappa_m)$ are found by formula (\ref{eq1201}), i.\ e.\
by solving a linear ODE with known coefficients and known initial conditions.

\item
Matrix ${\rm U}(k)$ is found at the points $k = \kappa_m$ by inverting formula (\ref{eq0809}).

\item
Functions $\hat U_0^j (k)$ are found for $k = \kappa_m$ by formula (\ref{u0tildej}).

\item
Functions $\hat U_0^j (k)$ are substituted into the embedding formula (\ref{embedding_a})  to get the function $\hat U_{0}(k,k_*)$.

\item
The directivity is found using formula (\ref{eq0617}) for the  points $\theta_m = \arccos(- \kappa_m / k_0)$.

\end{itemize}

One can see that all steps of this procedure can be done easily except the first one. To solve Problem~\ref{OE} we use
the technique introduced in \cite{Shanin2012}. Here we describe it.

Matrix $ \tilde {\rm M}_2(k)$ can be represented in the form:
\begin{equation}
{\rm \tilde M_2}(k) = {\rm H}(k)\left( \begin{array}{cc}
m(k - k_0) & 0 \\
0 & 1
\end{array} \right){\rm H}^{-1}(k),
\label{eigtildeM2}
\end{equation}
where $m(b)$, ${\rm H}(k)$ are introduced by (\ref{eq0901}) and (\ref{eq0908}), respectively.

\begin{figure}[ht]
\centerline{\epsfig{file=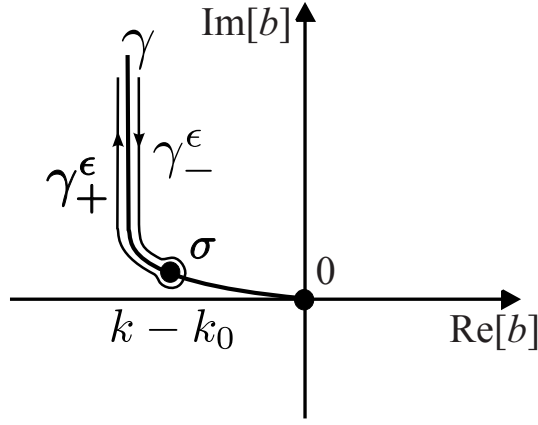}}
\caption{Contours $\gamma^{\epsilon}_{+},\sigma,\gamma^{\epsilon}_{-}$ }
\label{fig06b}
\end{figure}

The left--hand side of (\ref{eq1204}) can be rewritten as follows:
$$
{\rm OE}_{\gamma_+ \circ \gamma_-} \left[ \left(
\frac{{\rm r}(b)}{k - (k_0 + b)}
-
\frac{{\rm r^*}(b)}{k + (k_0 + b)}
\right) db \right]
$$
\begin{equation}
=
{\rm F}(k) {\rm OE}_{\sigma} \left[ \left(
\frac{{\rm r}(b)}{k - (k_0 + b)}
-
\frac{{\rm r^*}(b)}{k + (k_0 + b)}
\right) db \right]{\rm F}^{-1}(k),
\end{equation}
where $\sigma$ is a loop of a small radius $\epsilon$ encircling point $k - k_0$ in the positive direction, and
\begin{equation}
{\rm F}(k) = {\rm OE}_{\gamma_+^\epsilon} \left[ \left(
\frac{{\rm r}(b)}{k - (k_0 + b)}
-
\frac{{\rm r^*}(b)}{k + (k_0 + b)}
\right) db \right].
\end{equation}
Here $\gamma_+^\epsilon$ is contour going from $\i \infty$ to the start of the loop $\sigma$ along $\gamma_+$ (see Fig. \ref{fig06b}).
Let ${\rm r}(b)$ be represented in the form
\begin{equation}
{\rm r}(b) = {\rm P}(b)\left( \begin{array}{cc}
\xi_1(b) & 0 \\
0 & \xi_2(b)
\end{array} \right){\rm P}^{-1}(b).
\end{equation}
The columns of ${\rm P}(b)$ are the eigenvectors of $\rm r$. Almost everywhere matrix ${\rm P}(b)$ can be parametrized as follows:
\begin{equation}
{\rm P}(b)= \left( \begin{array}{cc}
1 & p_2(b) \\
p_1(b) & 1
\end{array} \right).
\end{equation}
One can see that as $\epsilon \to 0$
$$
{\rm OE}_{\sigma} \left[ \left(
\frac{{\rm r}(b)}{k - (k_0 + b)}
-
\frac{{\rm r^*}(b)}{k + (k_0 + b)}
\right) db \right] \to
$$
\begin{equation}
{\rm P}(k-k_0)\exp\left\{-2\pi\i\left( \begin{array}{cc}
\xi_1(k-k_0) & 0 \\
0 & \xi_2(k-k_0)
\end{array} \right)\right\}{\rm P}^{-1}(k-k_0).
\label{OErez}
\end{equation}
It follows from (\ref{OErez}) that the eigenvalues of ${\rm r}(b)$ are connected with the eigenvalues of $\tilde {\rm M}_2(k)$:
\begin{equation}
\label{ksi12}
\xi_1(k-k_0) = \frac{\i}{2\pi}\log (m(k-k_0)),\qquad
\xi_2(k) = 0.
\end{equation}
Thus, to find $\rm r$ one needs only to find $p_1(b), p_2(b)$.

\begin{figure}[ht]
\centerline{\epsfig{file=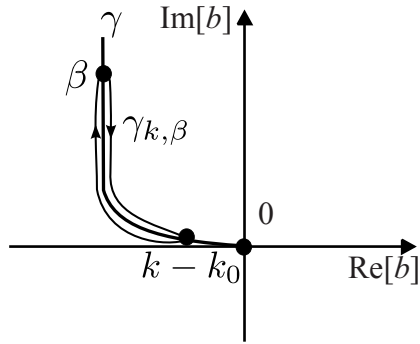}}
\caption{Contour $\gamma_{k,\beta}$ }
\label{fig06a}
\end{figure}

For $k,\beta \in \gamma, |\beta|> |k-k_0|$ define the function
\begin{equation}
{\rm K}(\beta,k) = {\rm OE}_{\gamma_{k,\beta}}\left[ \left(
\frac{{\rm r}(b)}{k - (k_0 + b)}
-
\frac{{\rm r^*}(b)}{k + (k_0 + b)}
\right) db \right]
\end{equation}
where contour $\gamma_{k,\beta}$ is shown in Fig.~\ref{fig06a}.
Function $  {\rm K} (\beta,k)$ obeys the equation
\begin{equation}
\frac{\ptl}{\ptl \beta}{\rm K}(\beta,k) = \frac{1}{k - (k_0 + \beta)}\left[{\rm r}(\beta),{\rm K}(\beta,k)\right] - \frac{1}{k + (k_0 + \beta)}\left[{\rm r}^*(\beta),{\rm K}(\beta,k)\right],
\label{VOE}
\end{equation}
where $[\cdot,\cdot]$ is a commutator.
Since ${\rm K}$ is adjoint to
\begin{equation}
{\rm OE}_{\sigma} \left[ \left(
\frac{{\rm r}(b)}{k - (k_0 + b)}
-
\frac{{\rm r^*}(b)}{k + (k_0 + b)}
\right) db \right],
\end{equation}
its eigenvalues are the same, it can be written as follows:
\begin{equation}
{\rm K}(\beta,k) = {\rm Q}(\beta,k)\left( \begin{array}{cc}
m(k-k_0) & 0 \\
0 & 1
\end{array} \right){\rm Q}^{-1}(\beta,k),
\end{equation}
where
\begin{equation}
{\rm Q}(\beta,k)= \left( \begin{array}{cc}
1 & q_2(\beta,k) \\
q_1(\beta,k) & 1
\end{array} \right).
\end{equation}
Taking $\beta \to k-k_0$ obtain the relation ${\rm Q}(k-k_0,k)= {\rm P}(k-k_0)$,
and thus
\begin{equation}
q_{1,2}(k-k_0,k)= p_{1,2}(k-k_0).
\label{riccond1}
\end{equation}

Taking $\beta \to \i\infty$ obtain the relation ${\rm Q}(i\infty,k)= {\rm H}(k)$, and thus
\begin{equation}
q_{1}(i\infty,k)= \alpha(k),\qquad q_2(i\infty,k) = 0,
\label{riccond2}
\end{equation}
where $\alpha$ is introduced by (\ref{eq0908}).

Elementary calculations demonstrate that equation (\ref{VOE}) is equivalent to the following system of two independent Riccati equations:
$$
\frac{\ptl q_{1,2}(\beta,k)}{\ptl \beta} = \frac{\xi_1(\beta)(p_{1,2}(\beta) - q_{1,2}(\beta,k))(1-p_{2,1}(\beta)q_{1,2}(\beta,k))}{(p_1(\beta)p_2(\beta)-1)(k - ( k_0 + \beta) )} +
$$
\begin{equation}
 \frac{\xi_1(\beta)(p_{2,1}(\beta) - q_{1,2}(\beta,k))(1-p_{1,2}(\beta)q_{1,2}(\beta,k))}{(p_1(\beta)p_2(\beta) - 1)(k + ( k_0 + \beta) )}.
\label{riccati}
\end{equation}

Thus we have to find $p_1(\beta), p_2(\beta)$
such that there exist solutions $q_{1,2}(\beta, k)$ of (\ref{riccati})
on the part of $\beta \in \gamma$ contained between the points $k- k_0$ and $i \infty$
obeying boundary conditions (\ref{riccond1}), (\ref{riccond2}).
This can be done using following numerical procedure. First, contour $\gamma$ should be meshed, i.\ e.\
an array of nodes $b_j, j =1 \dots N$ should be taken on it. Starting point $b_1$ should have the imaginary part big enough.  The end point $b_N$ is equal to $0$. At each point $b_j$ matrix ${\rm \tilde M_2}(b + k_0)$  is represented in the form (\ref{eigtildeM2}), i.\ e.\ the values $m(b)$ and $\alpha(b +k_0)$  are computed. The value $\xi_1(b + k_0)$ is computed by applying formula (\ref{ksi12}).

For the ``infinity'' point $b_1$ the following values are assigned:
\begin{equation}
p_{1}(b_1) = 0, \qquad  p_{2}(b_1) = 0.
\end{equation}
This is a natural choice for the asymptotics of the unknown coefficient, since ${\rm M}_2(k)$ tends to identity matrix, as
$k \to i\infty$.
The loop over $j=2...N$ is performed. At the $j$th step of the loop the values $p_{1,2}(b_j)$ are computed.
Thus on the $j$th step all values $p_{1,2}(b_1) \dots p_{1,2}(b_{j-1})$ are already found. On the $j$th step
of the loop equations (\ref{riccati}) are
solved on the contour $b \in (b_1,b_{j-1})$ for $q_{1,2}(b,b_j)$ using Runge--Kutta~4 method. The initial conditions are set as
\begin{equation}
q_{1}(b_1) = \alpha(b_1), \qquad  q_{2}(b_1) = 0
\end{equation}
following from (\ref{riccond2}). The values $q_{1,2}(b_1) \dots q_{1,2}(b_{j-1})$
are found. Then equations (\ref{riccati}) are solved on the segment $(b_{j-1},b_j)$ (only one step is performed) by Euler method. This method does not require the values of the
right--hand side at the end of the segment to be known. The values $q_{1,2}(b_j)$ are found. The assignment
\begin{equation}
p_{1,2}(b_j) = q_{1,2}(b_j)
\end{equation}
is performed, following from (\ref{riccond1}).

Thus matrix ${\rm {r}}(b)$ becomes known. It just remains to solve ODE1 (\ref{eq1001}) and then calculate antisymmetrical part of the directivity using the procedure described above. This is performed easily.

\begin{figure}[ht]
\centerline{\epsfig{file=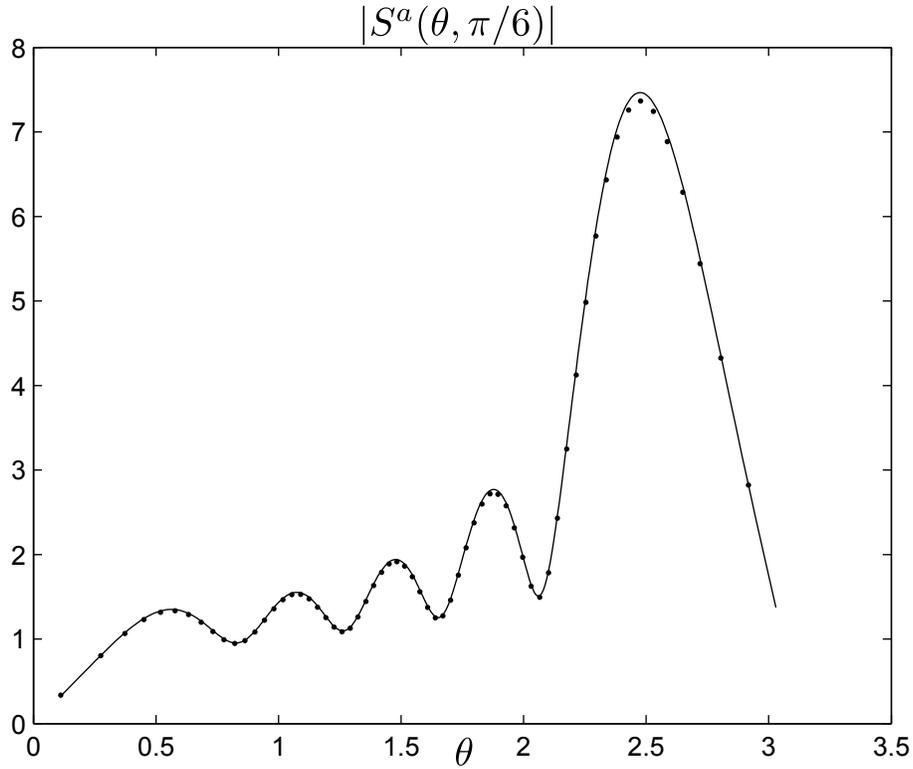}}
\caption{Dependence of $|S^{a}(\theta,\pi/6)|$ on $\theta$ for $ka = 8$ and $\eta = 1 - 0.25i$. Solid line corresponds to the method of integral equation and dotted line corresponds to the method of OE--equation}
\label{fig07}
\end{figure}

Numerical results are compared with solution obtained by the method of boundary integral equations (see Appendix~B). The dependence of $|S^{a}(\theta,\pi/6)|$ on $\theta$ for $ka = 8$ and $\eta = 1 - 0.25i$ is presented
in Fig.~\ref{fig07}. Solid line corresponds to the method of integral equation and dotted line corresponds to the
method of OE--equation. One can see that agreement is reasonable.


\subsection{Symmetrical case}
\begin{figure}[ht]
\centerline{\epsfig{file=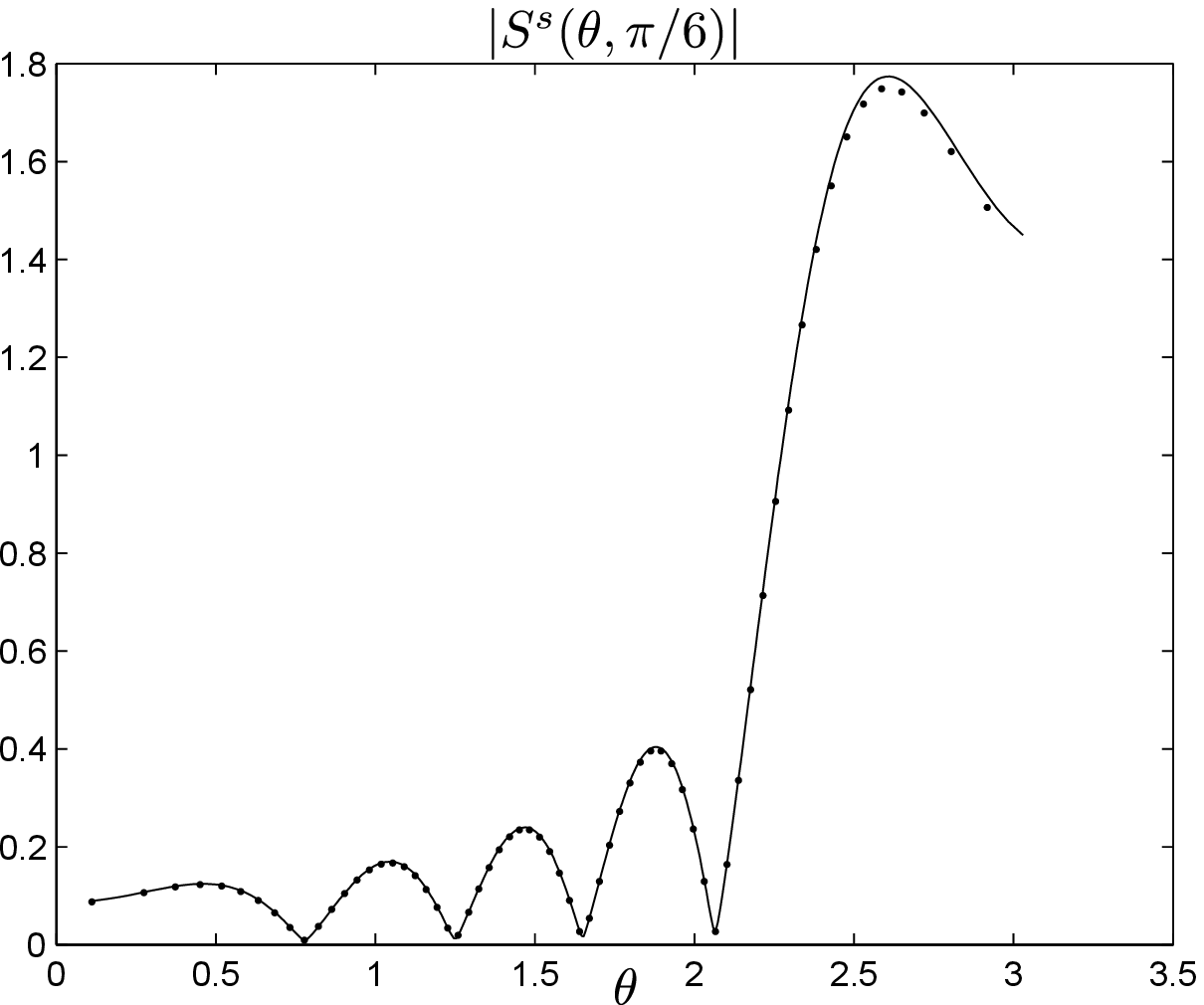}}
\caption{Dependence of $|S^{s}(\theta,\pi/6)|$ on $\theta$ for $ka = 8$ and $\eta = 1 - 0.25i$.
Solid line corresponds to the method of integral equation and dotted line corresponds to the method of OE--equation}
\label{fig08}
\end{figure}

The solution procedure in the symmetrical case is similar. Here we just present the final results. They are showed in
Fig.~\ref{fig08}. Dependence of $|S^{s}(\theta,\pi/6)|$ on $\theta$ for $ka = 8$ and $\eta = 1 - 0.25i$ is displayed.
Solid line corresponds to the method of integral equation and dotted line corresponds to the method of OE--equation.

\section{Conclusion}

In the current paper we present a new approach to matrix Riemann--Hilbert problems related to the
problem of diffraction by an impedance strip. The problems are of a quite general nature, so the
methods proposed here can potentially be applied to a wide class of problems.

The technique is based on an analytical result expressed in Theorem~1 and~2. The initial problem is
embedded into a family of similar problems indexed by parameter $b$, and it is shown that the
dependence of the solution on $b$ is described by an ordinary differential equation with a relatively
simple coefficient. Then, the Riemann--Hilbert problem is reformulated as a problem for an OE--equation,
i.\ e.\ a problem of reconstruction of the coefficients of an ODE by using the boundary data. There
is no analytical solution for the OE--problem in the general case, however some analytical technique
is available in the commutative case \cite{Shanin}. It is also worth to note that numerical solution
of the OE--problem can be very efficient since the problem is of Volterra nature (the unknown function on a contour
is found step by step).

To demonstrate the practical value of the analytical results obtained here we performed some computations of the
directivities for an impedance strip and compared the results with the integral equation method. The
agreement is nice, and this fact means for us mainly the validity of the method in general and the absence of
mistakes in the main formulae.
Here we do not pursue the aim to establish a robust and accurate numerical procedure based on the
new method.

\section*{Acknowledgements}

The work is supported by the grants
RFBR 14-02-00573, Scientific Schools 283.2014.2, and the Government grant 11.G34.31.0066.

The authors are grateful to participants of the seminar on wave diffraction held in S.Pb. branch of Steklov
Mathematical Institute of RAS (the chairman is Prof. V.\ M.\ Babich) for interesting discussions.

\section*{{\bf Appendix A. Index of Riemann - Hilbert problem}}
Let us prove formula (\ref{eq0903}). Obviously,
\begin{equation}
{\rm Idx} = \i {\rm Arg}[m(b)]|^{\i \infty}_{0},
\end{equation}
where $ {\rm Arg}[\cdot]$ is the argument of the function.
Introduce
\begin{equation}
f_1(b) = \i \sqrt{k_0^2 - (k_0 + b)^2} + \eta,
\end{equation}
\begin{equation}
f_2(b) = \i \sqrt{k_0^2 - (k_0 + b)^2} - \eta.
\end{equation}
Thus
\begin{equation}
{\rm Idx} = \i {\rm Arg}[f_1(b)]|^{\i \infty}_{0} - \i {\rm Arg}[f_2(b)]|^{\i \infty}_{0}.
\end{equation}
Here we consider only the case ${\rm Re}[\eta] > 0$ for which no deformation of contour $\gamma$ is needed,
Due to the rules of analytical continuation our proof is valid for all $\eta$ lying in the lower half--plane.
 \begin{figure}[ht]
\centerline{\epsfig{file=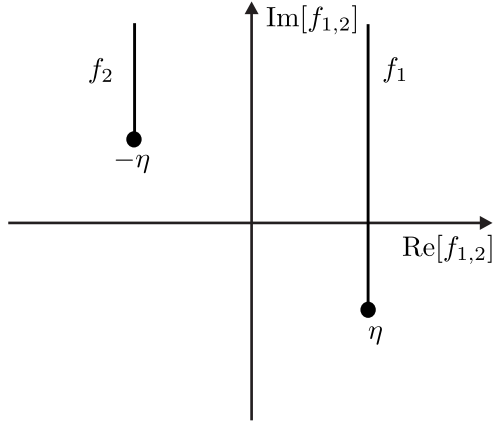}}
\caption{Contours $f_{1,2}(b)$ for $b \in \gamma$ }
\label{fig09}
\end{figure}

One can notice(see Fig.~\ref{fig09}) that the argument of $f_{1,2}(b)$
changes from ${\rm Arg}[\pm \eta]$ to $\pi/2$ while $b$ goes from $0$ to $\i \infty$ along  $\gamma$. Therefore
\begin{equation}
{\rm Idx} = \i(\pi/2 - {\rm Arg}[\eta]) - i(\pi/2 - {\rm Arg}[-\eta]) =
i({\rm Arg}[-\eta] - Arg[\eta]) = \i \pi.
\end{equation}


\section*{{\bf Appendix B. Integral equation method}}

{\bf The antisymmetrical case.}
To verify the results  obtained by  the OE--equation method we also solved the problem of diffraction by an impedance strip using the integral equation method. In the antisymmetrical case one can obtain following equation with the help of double layer potential:
\begin{equation}
\label{int_eq_a}
\left(\frac{\ptl^2}{\ptl x^2} + k_0^2\right)\int^{a}_{-a}G(x-x',0)\nu(x')dx' + \frac{1}{2}\eta \nu(x) = \i k_0 \sin\theta^{in}\exp\{-\i k_0 x \cos \theta^{in}\},
\end{equation}
where
$$
G(x,y) = -\frac{\i}{4}H_0^{(1)}(k_0\sqrt{x^2 + y^2}),
$$
$H_0^{(1)}(z)$ is a Hankel function of the first kind, $\nu$ is a double--layer potential:
\begin{equation}
u^a(x,y) = -\int^a_{-a}\frac{\ptl}{\ptl y}G(x-x',y)\nu(x')dx',
\end{equation}
$u^{\rm a}$ is the antisymmetrical part of the scattered field $u^{\rm sc}$.
On the strip, $\nu$  is connected with $u^{\rm a}(x,y)$ by a simple formula:
$$
u^{\rm a}(x,+0) = -\frac{1}{2}\nu(x), \quad -a<x<a.
$$
Antisymmetrical part of the directivity can be calculated using formula
\begin{equation}
S^{\rm a} (\theta,\theta^{\rm in}) = -e^{- i \pi /4} k_0 \sin \theta \, \int \limits_{-a}^a u^{\rm a}(x,+0) e^{-\i k_0 x \cos\theta } dx.
\end{equation}

Problem (\ref{int_eq_a}) can be discretized and solved numerically with help of the standard techniques.


{\bf Symmetrical case.}
In the symmetrical case it is natural to use a single layer potential. One can obtain the  following
integral equation:
\begin{equation}
\label{int_eq_s}
 \frac{1}{2} \mu(x) - \eta\int^{a}_{-a}G(x-x',0)\mu(x')dx' =\eta \exp\{-\i k_0 x\cos \theta^{in}\},
\end{equation}
where $\mu$ is a single layer potential:
\begin{equation}
u^s(x,y) = \int^a_{-a}G(x-x',y)\mu(x')dx'.
\end{equation}
Here $u^{\rm s}$ is the symmetrical part of the scattered field $u^{\rm sc}$.
The normal derivative of the field
on the strip is connected with $\mu(x)$ as follows:
 $$
\frac{\ptl u^{\rm s}}{\ptl y}(x,+0) = \frac{1}{2}\mu(x), \quad -a<x<a.
$$
Symmetrical part of the directivity can be calculated using the formula:
\begin{equation}
S^{\rm s} (\theta,\theta^{\rm in}) = e^{- i \pi /4} \, \int \limits_{-a}^a \frac{\ptl u^{\rm s}}{\ptl y}(x,+0) e^{-\i k_0 x \cos\theta } dx.
\end{equation}

\end{document}